\documentclass[11pt]{article}
\usepackage{latexsym,amssymb, color}
\usepackage{psfrag}
\usepackage{amsmath}

\newtheorem{Theorem}{\bf Theorem}[section]
\newtheorem{Lemma}{\bf Lemma}[section]
\newtheorem{Proposition}{\bf Proposition}[section]  
\newtheorem{Corollary}{\bf Corollary}[section]
\newtheorem{Remark}{\bf Remark}[section]
\newtheorem{Example}{\bf Example}[section]
\newtheorem{Definition}{\bf Definition}[section]

\newenvironment{theorem}{\begin{Theorem}$\!\!\!$}{\end{Theorem}}
\newenvironment{lemma}{\begin{Lemma}$\!\!\!$}{\end{Lemma}}
\newenvironment{proposition}{\begin{Proposition}$\!\!\!$}{\end{Proposition}}
\newenvironment{corollary}{\begin{Corollary}$\!\!\!$}{\end{Corollary}}
\newenvironment{remark}{\begin{Remark}$\!\!\!$}{\end{Remark}}

\newenvironment{definition}{\begin{Definition}$\!\!\!$}{\end{Definition}}

\numberwithin{equation}{section}

\begin{document}
\title{Hot spots of solutions to the heat equation\\
with inverse square potential}
\author{Kazuhiro Ishige, 
Yoshitsugu Kabeya and Asato Mukai
}
\date{}
\maketitle
\begin{abstract}
We investigate the large time behavior of the hot spots of the solution 
to the Cauchy problem 
$$
\left\{
\begin{array}{ll}
\partial_t u-\Delta u+V(|x|)u=0 &\mbox{in}\quad{\bf R}^N\times(0,\infty),\vspace{3pt}\\
u(x,0)=\varphi(x) & \mbox{in}\quad{\bf R}^N,
\end{array}
\right.
$$
where $\varphi\in L^2({\bf R}^N,e^{|x|^2/4}dx)$ and $V=V(r)$ decays quadratically as $r\to\infty$. 
In this paper, based on the arguments in [K. Ishige and A. Mukai, preprint (arXiv:1709.00809)], 
we classify the large time behavior of the hot spots of $u$ and  reveal 
the relationship between the behavior of the hot spots and the harmonic functions for $-\Delta+V$. 
\end{abstract}
\vspace{40pt}
\noindent Addresses:

\smallskip
\noindent K. I.:  Mathematical Institute, Tohoku University,
Aoba, Sendai 980-8578, Japan.\\
\noindent 
E-mail: {\tt ishige@m.tohoku.ac.jp}\\

\smallskip
\noindent 
Y. K.: Department of Mathematical Sciences, Osaka Prefecture University, 
Sakai 599-8531, Japan. \\
\noindent 
E-mail: {\tt kabeya@ms.osakafu-u.ac.jp}\\

\smallskip
\noindent 
A. M.: Mathematical Institute, Tohoku University,
Aoba, Sendai 980-8578, Japan.\\
\noindent 
E-mail: {\tt asato.mukai.t7@dc.tohoku.ac.jp}\\
\newpage
\section{Introduction}
Let $u$ be a solution of 
\begin{equation}
\label{eq:1.1}
\left\{
\begin{array}{ll}
\partial_t u-\Delta u+V(|x|)u=0 & \mbox{in}\quad{\bf R}^N\times(0,\infty),\vspace{3pt}\\
u(x,0)=\varphi(x) & \mbox{in}\quad{\bf R}^N,
\end{array}
\right.
\end{equation}
where $\varphi\in L^2({\bf R}^N,e^{|x|^2/4}\,dx)$. 
Here $L_V:=-\Delta+V$ is a nonnegative Schr\"odinger operator on $L^2({\bf R}^N)$, 
where $N\ge 2$ and $V$ is a radially symmetric inverse square potential. 
More precisely, we assume the following condition (V) on the potential $V$:  
\begin{equation*}
(V)\qquad
\left\{
\begin{array}{ll}
({\rm i}) & \mbox{$V=V(r)\in C^1((0,\infty))$};\vspace{7pt}\\
({\rm ii}) &  \displaystyle{\lim_{r\to 0}r^{-\theta}\left|r^2V(r)-\lambda_1\right|=0},
 \quad
 \displaystyle{\lim_{r\to\infty}r^\theta\left|r^2V(r)-\lambda_2\right|=0},\vspace{3pt}\\
  & \mbox{for some $\lambda_1$, $\lambda_2\in[\lambda_*,\infty)$ with $\lambda_*:=-(N-2)^2/4$ and $\theta>0$};\vspace{7pt}\\
({\rm iii}) &  \displaystyle{\sup_{r\ge 1}\left|r^3V'(r)\right|<\infty}. 
\end{array}
\right.
\qquad
\end{equation*}
Nonnegative Schr\"odinger operators and their heat semigroups appear in various fields 
and have been studied intensively 
by many authors since the pioneering work due to Simon~\cite{S} 
(see e.g., \cite{CK1}--\cite{Gri}, \cite{IK02}--\cite{IM}, \cite{LS}--\cite{PZ}, 
\cite{VZ}--\cite{Zhang} and references therein).  
The inverse square potential is a typical one appearing in the study of the Schr\"odinger operators 
and it also arises in the linearized analysis for nonlinear diffusion equations, 
in particular, in solid-fuel ignition phenomena which can be modeled by 
$\partial_t u=\Delta u+e^u$. 
(See e.g., \cite{BE}, \cite{Marchi}, \cite{V}, \cite{VZ} and \cite{Z}.) 
In this paper we investigate the large time behavior of the hot spots 
$$
H(u(t)):=\biggr\{x\in{\bf R}^N\,:\,u(x,t)=\sup_{y\in{\bf R}^N}u(y,t)\biggr\}
$$
for the solution~$u$ of \eqref{eq:1.1}. 
The study of the large time behavior of the hot spots is delicate and 
it is obtained by the higher order asymptotic expansion of the solutions. 
\vspace{3pt}

We say that $L_V:=-\Delta+V(|x|)$ is nonnegative on $L^2({\bf R}^N)$ if 
$$
\int_{{\bf R}^N}\left[|\nabla\phi|^2+V(|x|)\phi^2\right]\,dx\ge 0,
\qquad\phi\in C_0^\infty({\bf R}^N\setminus\{0\}). 
$$
When $L_V$ is nonnegative,  
we say that 
\begin{itemize}
  \item 
  $L_V$ is subcritical if, for any $W\in C_0^\infty({\bf R}^N)\setminus\{0\}$, 
  $L-\epsilon W$ is nonnegative for all sufficiently small $\epsilon>0$; 
  \item 
  $L_V$ is critical if $L_V$ is not subcritical. 
\end{itemize}
For $\lambda\ge\lambda_*$, let $A^\pm(\lambda)$ be the roots of the algebraic equation $\alpha^2+(N-2)\alpha-\lambda=0$ 
such that $A^-(\lambda)\le A^+(\lambda)$, that is 
\begin{equation}
\label{eq:1.2}
A^\pm(\lambda):=\frac{-(N-2)\pm\sqrt{(N-2)^2+4\lambda}}{2}. 
\end{equation}
Under condition~(V), 
there exists a unique solution $U$ of 
\begin{equation}
\label{eq:1.3}
U''+\frac{N-1}{r}U'-V(r)U=0\quad\mbox{in}\quad(0,\infty)
\quad\mbox{with}\quad
\lim_{r\to 0}r^{-A^+(\lambda_1)}U(r)=1. 
\end{equation}
Assume that $L_V$ satisfies one of the following three conditions :
$$
\begin{array}{ll}
& {\rm (S)}\,:\mbox{$L$ is subcritical and $\lambda_2>\lambda_*$};
\qquad\qquad
{\rm (S_*)}:\mbox{$L$ is subcritical and $\lambda_2=\lambda_*$};\vspace{3pt}\\
 &{\rm (C)}:\mbox{$L$ is critical and $A^-(\lambda_2)>-N/2$}.
\end{array}
$$ 
Then $U>0$ in $(0,\infty)$
and $U$ satisfies 
\begin{equation}
\label{eq:1.4}
U(r)\thicksim c_*v(r)\quad\mbox{as}\quad r\to\infty
\end{equation}
for some positive constant $c_*$, where 
$$
v(r):=
\left\{
\begin{array}{ll}
r^A & \mbox{in the cases of $\rm (S)$ and $\rm (C)$},\vspace{3pt}\\
r^A\log(2+r) & \mbox{in the case of $\rm (S_*)$}.
\end{array}
\right.
$$
Here $A=A^+(\lambda_2)$ in the cases of (S) and $\rm (S_*)$ 
and $A=A^-(\lambda_2)$ in the case of (C). 
(See \cite[Theorem 1.1]{IKO}.) 
Furthermore, the following result holds. 
See \cite[Theorem~1.4]{IM}.
\begin{theorem}
\label{Theorem:1.1}
Let $N\ge 2$ and $\varphi\in L^2({\bf R}^N,e^{|x|^2/4}\,dx)$. 
Let $u$ be a solution of \eqref{eq:1.1} under condition~{\rm (V)} with $A>-N/2$. 
Set
$$
M(\varphi):=\frac{1}{c_*\kappa}\int_{{\bf R}^N}\varphi(x)U(|x|)\,dx,
\quad
\kappa:=2^{N+2A}\pi^{\frac{N}{2}}\Gamma\left(\frac{N+2A}{2}\right)\biggr/\Gamma\left(\frac{N}{2}\right),
$$
where $\Gamma$ is the Gamma function. 
\begin{itemize}
  \item[{\rm (a)}] 
  In cases {\rm (S)} and {\rm (C)}, 
  $$
  \lim_{t\to\infty}t^{\frac{N+A}{2}}u(t^{\frac{1}{2}}y,t)=M(\varphi)|y|^Ae^{-\frac{|y|^2}{4}}
  $$
  in $L^2({\bf R}^N,e^{|y|^2/4}\,dy)$ and in $L^\infty(K)$ for any compact set $K\subset{\bf R}^N\setminus\{0\}$. 
  Furthermore, for any sufficiently small $\epsilon>0$,
  $$
  t^{\frac{N+2A}{2}}\frac{u(x,t)}{U(|x|)}=c_*^{-1}M(\varphi)+o(1)
  +t^{-1}O(|x|^2)\quad\mbox{as}\quad t \to \infty
  $$
  uniformly for $x\in{\bf R}^N$ with $|x|\le\epsilon(1+t)^{\frac{1}{2}}$. 
  \item[{\rm (b)}] 
  In case {\rm ($\mbox{S}_*$)},
  $$
  \lim_{t\to\infty}t^{\frac{N+A}{2}}(\log t)u(t^{\frac{1}{2}}y,t)=2M(\varphi)|y|^Ae^{-\frac{|y|^2}{4}}
  $$
  in $L^2({\bf R}^N,e^{|y|^2/4}\,dy)$ and in $L^\infty(K)$ for any compact set $K\subset{\bf R}^N\setminus\{0\}$. 
  Furthermore, for any sufficiently small $\epsilon>0$,
  $$
  t^{\frac{N+2A}{2}}(\log t)^2\,\frac{u(x,t)}{U(|x|)}=4c_*^{-1}M(\varphi)+o(1)+O(t^{-1}|x|^2) \quad\mbox{as}\quad t \to \infty
  $$
  uniformly for $x\in{\bf R}^N$ with $|x|\le\epsilon(1+t)^{\frac{1}{2}}$. 
\end{itemize}
\end{theorem}
(See Section~2.2.)
In this paper, combing the arguments in \cite{IK02}--\cite{IK05} and \cite{IM}, 
we study the large time behavior of $H(u(t))$ in the cases (S), (${\rm S_*}$) and (C) 
and reveal the relationship between the large time behavior of $H(u(t))$ and the corresponding harmonic functions. 
We remark that $L_V$ is not necessarily subcritical. 
\vspace{3pt}

The behavior of the hot spots for parabolic equations 
in unbounded domains has been studied since the pioneering work by Chavel and Karp~\cite{CK}, 
who studied the behavior of the hot spots for the heat equations on some non-compact Riemannian manifolds. 
In particular, for the heat equation on ${\bf R}^N$ with nonnegative initial data $\varphi\in L^\infty_c({\bf R}^N)$, 
they proved: 
\begin{itemize}
  \item[(H1)] $H(e^{t \Delta}\varphi)$ is a subset of the closed convex hull of the support of the initial function $\varphi$;
  \item[(H2)] There exists $T>0$ such that $H(e^{t \Delta}\varphi)$ consists of only one point and moves along a smooth curve for any $t\ge T$;
  \item[(H3)] $\displaystyle{\lim_{t\to\infty}H(e^{t \Delta}\varphi)=\int_{{\bf R}^N}\,x\varphi(x)\,dx\biggr/\int_{{\bf R}^N}\varphi(x)\,dx}$. 
\end{itemize}
(See also Remark~\ref{Remark:3.3}.)
The behavior of the hot spots for the heat equation on the half space of ${\bf R}^N$ and on the exterior domain of a ball 
was studied in \cite{I1}, \cite{I2} and \cite{JS}. 
Subsequently, in~\cite{IK02}--\cite{IK05}, 
the first and the second authors of this paper developed the arguments in \cite{I1} and \cite{I2} 
and studied the large time behavior of the hot spots for the solution of \eqref{eq:1.1} under condition~(V) 
in the subcritical case with some additional assumptions. 
\vspace{3pt}

Our arguments in this paper are based on \cite{IM}, 
where the precise description of the large time behavior of the solution of \eqref{eq:1.1} was discussed under condition~(V). 
Applying the arguments in \cite{IM}, we modify the arguments in \cite{IK02}--\cite{IK05} 
and study the large time behavior of the hot spots. 
We study the following subjects
when the hots spots tend to the space infinity as $t\to\infty$:
\begin{itemize}
  \item[(a)] 
  The rate and the direction for the hot spots to tend to the space infinity as $t\to\infty$;
  \item[(b)] 
  The number of  the hot spots for sufficiently large $t$.
\end{itemize}
On the other hand, when the hots spots accumulate to a point $x_*$, 
we characterize the limit point $x_*$ by the positive harmonic function~$U$. 
Furthermore, we give a sufficient condition for the hot spots to 
consist of only one point and to move along a smooth curve. 
\vspace{3pt}

The rest of this paper is organized as follows.  
In Section~2 we formulate a definition of the solution of \eqref{eq:1.1}. 
Furthermore, we recall some preliminary results on the behavior of the solution of \eqref{eq:1.1} 
and prove some lemmas. 
In Section~3 we study the large time behavior of the hot spots for problem~\eqref{eq:1.1}. 
\section{Preliminaries}
In this section we formulate the definition of the solution of \eqref{eq:1.1} 
and recall some results on the behavior of the solution. 
Throughout this paper, 
for positive functions $f=f(s)$ and $g=g(s)$ in $(R,\infty)$ for some $R>0$, 
we say that $f(s)\thicksim g(s)$ for all sufficiently large $s>0$ if $\lim_{s\to\infty}f(s)/g(s)=1$. 
Furthermore, we say that $f(s)\asymp g(s)$ for all sufficiently large $s>0$ if 
if there exists $C>0$ such that $C^{-1}\le f(s)/g(s)\le C$ for all sufficiently large $s>0$.  
\vspace{3pt}

Assume condition~(V) and let $U$ be a positive solution of \eqref{eq:1.3}. 
It follows from \eqref{eq:1.2} and \eqref{eq:1.3} that $U^2\in L^1_{\rm loc}({\bf R}^N)$.
Consider the Cauchy problem 
\begin{equation}
\tag{P}
\qquad
\left\{
\begin{array}{ll}
\partial_t u_*+L_*u_*=0 & \quad\mbox{in}\quad{\bf R}^N\times(0,\infty),
\vspace{5pt}\\
u_*(x,0)=\varphi_*(x) & \quad\mbox{in}\quad{\bf R}^N,
\vspace{3pt}
\end{array}
\right.
\end{equation}
where 
$$
L_*u_*:=-\frac{1}{\nu}\mbox{div}\,(\nu\nabla u_*),
\quad
\nu:=U^2\in L^1_{\rm loc}({\bf R}^N),
\quad
\varphi_*\in L^2({\bf R}^N,\,\nu\,dx).
$$
\begin{definition}
\label{Definition:2.1}
Let $\varphi_*\in L^2({\bf R}^N,\,\nu\,dx)$. 
We say that $u_*$ is a solution of $(P)$ if 
\begin{equation*}
\begin{split}
  & u_*\in C([0,\infty):L^2({\bf R}^N,\,\nu\,dx))\,\cap\,L^2((0,\infty):H^1({\bf R}^N,\,\nu\,dx)),\vspace{5pt}\\
  & \int_0^\infty\int_{{\bf R}^N}\left[-u_*\partial_th+\nabla u_*\nabla h\right]\nu\,dx\,d\tau=0
\quad\mbox{for any $h\in C^\infty_0({\bf R}^N\times(0,\infty))$},\\
  & \lim_{t\to +0}\|u_*(t)-\varphi_*\|_{L^2({\bf R}^N,\,\nu\,dx)}=0.
\end{split}
\end{equation*}
\end{definition}
Problem~(P) possesses a unique solution $u_*$ such that 
\begin{equation}
\label{eq:2.1}
\|u_*(t)\|_{L^2({\bf R}^N,\,\nu\,dx)}\le\|\varphi_*\|_{L^2({\bf R}^N,\,\nu\,dx)},
\qquad
u_*(\cdot,t)\in C({\bf R}^N),
\end{equation}
for any $t>0$. See \cite[Section~2]{IM}. 
We state the definition of the solution of \eqref{eq:1.1}. 
\begin{definition}
\label{Definition:2.2}
Let $u$ be a measurable function in ${\bf R}^N\times(0,\infty)$ and $\varphi\in L^2({\bf R}^N)$. 
Define $u_*(x,t):=u(x,t)/U(|x|)$ and $\varphi_*(x):=\varphi(x)/U(|x|)$. 
Then we say that $u$ is a solution of \eqref{eq:1.1} if $u_*$ is a solution of {\rm (P)}. 
\end{definition}
\subsection{Asymptotic behavior of solutions}
In this subsection we recall some results 
on the large time behavior of radially symmetric solutions of \eqref{eq:1.1}. 
See \cite[Theorems~1.1, 1.2 and 1.3]{IM}.
\begin{proposition}
\label{Proposition:2.1}
Let $N\ge 2$ and assume condition~$(V)$. 
Let $L_V$ satisfy either {\rm (S)} or {\rm (C)}. 
Let $u=u(|x|,t)$ be a radially symmetric solution of \eqref{eq:1.1} 
such that $\varphi\in L^2({\bf R}^N,e^{|x|^2/4}\,dx)$. 
\begin{itemize}
  \item[{\rm (a)}]
  Define $w=w(\xi,s)$ by 
  $$
  w(\xi,s):=(1+t)^{\frac{N+A}{2}}u(r,t)\quad\mbox{with}\quad \xi=(1+t)^{-\frac{1}{2}}r\ge 0,\,\,\, s=\log(1+t)\ge 0.
  $$
  Then there exists a positive constant $C$ such that 
  $$
  \sup_{s>0}\|w(s)\|_{L^2({\bf R}^N,e^{|\xi|^2/4}\,d\xi)}\le C\|w(0)\|_{L^2({\bf R}^N,e^{|\xi|^2/4}\,dx)}.
  $$
  Furthermore, 
  $$
  \lim_{s\to\infty}w(\xi,s)=c_Am(\varphi)|\xi|^Ae^{-\frac{|\xi|^2}{4}}
  \quad\mbox{in}\quad L^2({\bf R}^N,e^{|\xi|^2/4}\,dx)\,\cap\,C^2(K) 
  $$
  for any compact set $K$ in ${\bf R}^N\setminus\{0\}$, 
  where
  $$
  c_A=\biggr[2^{N+2A-1}\Gamma\left(\frac{N+2A}{2}\right)\biggr]^{-1/2},
  \quad
  m(\varphi):=\frac{c_A}{c_*}\int_0^\infty \varphi(r)U(r)r^{N-1}\,dr.
  $$
  Moreover, if $m(\varphi)=0$, then 
  $$
  \|w(s)\|_{L^2({\bf R}^N,e^{|\xi|^2/4}\,dx)}+\|w(s)\|_{C^2(K)}=O(e^{-s})
  \quad\mbox{as}\quad s\to\infty. 
  $$
  \item[{\rm (b)}]  
  Set $u_*(r,t):=u(r,t)/U(r)$. 
  Then, for any $j\in\{0,1,2\dots\}$, $\partial_t^j u_*\in C([0,\infty)\times(0,\infty))$. 
  Define 
  $$
  G(r,t):=u_*(r,t)-\left[u_*(0,t)+(\partial_tu_*)(0,t)F(r)\right]\quad\mbox{for}\quad r\in[0,\infty),\,\,t>0,
  $$
  with
  \begin{equation}
  \label{eq:2.2}
  F(r):=\int_0^r s^{1-N}[U(s)]^{-2}\left(\int_0^s \tau^{N-1}U(\tau)^2\,d\tau\right)\,ds.
  \end{equation}
  Then 
  $$
  \lim_{t\to\infty}\,t^{\frac{N+2A}{2}}u_*(0,t)=\frac{c_A}{c_*}m(\varphi),\qquad
  \lim_{t\to\infty}t^{\frac{N+2A}{2}+1}(\partial_tu_*)(0,t)=-\frac{(N+2A)c_A}{2c_*}m(\varphi).
  $$
  Furthermore, for any $T>0$ and any sufficiently small $\epsilon>0$, 
  there exists $C_1>0$ such that 
  $$
  |(\partial_r^\ell G)(r,t)|\le C_1t^{-\frac{N+2A}{2}-2}r^{4-\ell}\|\varphi\|_{L^2({\bf R}^N,e^{|x|^2/4}\,dx)}
  $$
  for $\ell\in\{0,1,2\}$, $0\le r\le\epsilon(1+t)^{\frac{1}{2}}$ and $t\ge T$.  
\end{itemize}
\end{proposition}
\begin{proposition}
\label{Proposition:2.2}
Let $N\ge 2$ and assume condition~$(V)$. 
Let $L_V$ satisfy {\rm ($\mbox{S}_*$)}. 
Let $u=u(|x|,t)$ be a radially symmetric solution of \eqref{eq:1.1} 
such that $\varphi\in L^2({\bf R}^N,e^{|x|^2/4}\,dx)$. 
Then $d=N+2A=2$ and the following holds. 
\begin{itemize}
  \item[{\rm (a)}] 
  Let $w$ be as in Proposition~{\rm\ref{Proposition:2.1}} and $K$ a compact set in ${\bf R}^N\setminus\{0\}$. 
  Then there exists a positive constant $C_1$ such that 
  $$
  \sup_{s>0}\,(1+s)\|w(s)\|_{L^2({\bf R}^N,e^{|\xi|^2/4}\,dx)}\le C_1\|w(0)\|_{L^2({\bf R}^N,e^{|\xi|^2/4}\,dx)}.
  $$
  Furthermore, 
  $$
  \lim_{s\to\infty}sw(\xi,s)=2c_Am(\varphi)|\xi|^Ae^{-\frac{|\xi|^2}{4}}
  \quad\mbox{in}\quad
  L^2({\bf R}^N,e^{|\xi|^2/4}\,dx)\,\cap\,C^2(K),
  $$
  where $c_A$ and $m(\varphi)$ are as in Proposition~{\rm\ref{Proposition:2.1}}.
  \item[{\rm (b)}] 
  Let $u_*$ be as in Proposition~{\rm\ref{Proposition:2.1}}. 
  Then 
  \begin{equation*}
  \begin{split}
   & \lim_{t\to\infty}t(\log t)^2 u_*(0,t)=2\sqrt{2}c_*^{-1}m(\varphi),\\
   & \lim_{t\to\infty}t^2(\log t)^2(\partial_t u_*)(0,t)=-2\sqrt{2}c_*^{-1}m(\varphi).
  \end{split}
  \end{equation*}
  Furthermore, for any $T>0$ and any sufficiently small $\epsilon>0$, 
  there exists $C_2>0$ such that 
  $$
  |(\partial_r^\ell G_2)(r,t)|\le C_2t^{-3}[\log(2+t)]^{-2}r^{4-\ell}\|\varphi\|_{L^2({\bf R}^N,e^{|x|^2/4}\,dx)}
  $$
  for $\ell\in\{0,1,2\}$, $0\le r\le\epsilon(1+t)^{\frac{1}{2}}$ and $t\ge T$. 
  \end{itemize}
\end{proposition}

Following \cite[Section~1.2]{IM}, 
we apply Propositions~\ref{Proposition:2.1} and \ref{Proposition:2.2} to 
obtain the large time behavior of the solution of \eqref{eq:1.1}. 
Let $\{\omega_k\}_{k=0}^\infty$ be 
the eigenvalues of 
$$
-\Delta_{{\bf S}^{N-1}}Q=\omega Q\quad\mbox{on}\quad{\bf S}^{N-1},
\qquad
Q\in L^2({\bf S}^{N-1}), 
$$
where $\Delta_{{\bf S}^{N-1}}$ be the Laplace-Beltrami operator on the unit sphere ${\bf S}^{N-1}$. 
Then $\omega_k:=k(N+k-2)$ for $k=0,1,2,\dots$. 
Let
$\ell_k$  and $\{Q_{k,i}\}_{i=1}^{\ell_k}$ be 
the dimension and the orthonormal system of the eigenspace corresponding to $\omega_k$, respectively. 
In particular, $\ell_0=1$, $\ell_1=N$ and 
\begin{equation}
\label{eq:2.3}
Q_{0,1}\left(\frac{x}{|x|}\right)=q_*:=|{\bf S}^{N-1}|^{-\frac{1}{2}},\qquad
Q_{1,i}\left(\frac{x}{|x|}\right)=q_N\frac{x_i}{|x|}
\quad\mbox{with}\quad
q_N=N^{\frac{1}{2}}q_*,
\end{equation}
where $i=1,\dots,N$. Here $|{\bf S}^{N-1}|$ is the surface area of ${\bf S}^{N-1}$, 
that is $|{\bf S}^{N-1}| = 2\pi^{N/2} / \Gamma(N/2).$
For any $\varphi\in L^2({\bf R}^N,e^{|x|^2/4}\,dx)$, 
we can find radially symmetric functions $\{\phi^{k,i}\}\subset L^2({\bf R}^N,e^{|x|^2/4}\,dx)$ such that 
$$
\varphi=\sum_{k=0}^\infty\sum_{i=1}^{\ell_k} \varphi^{k,i}\quad\mbox{in}\quad L^2({\bf R}^N,e^{|x|^2/4}\,dx),
\quad
\varphi^{k,i}(x):=\phi^{k,i}(|x|)Q_{k,i}\left(\frac{x}{|x|}\right)
$$
(see \cite{I1} and \cite{IK02}). 
Let $u(x,t):=[e^{-tL_V}\varphi](x)$ and $u_0(x,t):=[e^{-tL_V}\varphi^{0,1}](|x|)$. 
For $m=1,2,3,\dots$, define
$$
R_m(x,0):=\varphi(x)-\sum_{k=0}^{m-1}\sum_{i=1}^{\ell_k}\varphi^{k,i}(x),\qquad
R_m(x,t):=[e^{-tL_V}R_m(0)](x).
$$
Then we have
\begin{equation}
\label{eq:2.4}
\begin{split}
 & u_{k,i}(x,t):=[e^{-tL_V}\varphi^{1,i}](|x|)=[e^{-tL_V^k}\phi^{k,i}](|x|)\,Q_{k,i}\left(\frac{x}{|x|}\right),\\
 & u(x,t)=u_0(x,t)+\sum_{k=1}^{m-1}\sum_{i=1}^{\ell_k}u_{k,i}(x,t)+R_m(x,t),
\end{split}
\end{equation}
for $x\in{\bf R}^N$ and $t>0$, where 
$L_V^k:=-\Delta+V_k(|x|)$ and $V_k(|x|):=V(|x|)+\omega_k|x|^{-2}$. 
Here $V_k$ satisfies condition~(V) with $\lambda_1$ and $\lambda_2$ 
replaced by $\lambda_1+\omega_k$ and $\lambda_2+\omega_k$, respectively, 
and the operator $L_V^k$ is subcritical if $k\ge 1$.
Let $U_k$ be a (unique) solution of \eqref{eq:1.3} with $V$ replaced by $V_k$.
Then 
\begin{equation}
\label{eq:2.5}
\begin{split}
 & U_k(r)\thicksim r^{A^+(\lambda_1+\omega_k)}\quad\mbox{as}\quad r\to+0,\\
 & U_k(r)\thicksim c_k\,r^{A^+(\lambda_2+\omega_k)},
 \quad U_k'(r)=O\left(r^{A^+(\lambda_2+\omega_k)-1}\right)\quad\mbox{as}\quad r\to\infty,
\end{split}
\end{equation}
for some positive constant $c_k$. 
Then, by Propositions~\ref{Proposition:2.1} and \ref{Proposition:2.2} 
we obtain the precise description 
of the large time behavior of $e^{-tL_V^k}\phi^{k,i}$, where $k=0,1,2,\dots$. 
In particular, for $k=0,1,2,\dots$, 
for any sufficiently small $\epsilon>0$, 
we have
\begin{equation}
\label{eq:2.6}
\begin{split}
 & t^{\frac{N+2A_k}{2}}\partial_r^\ell\frac{[e^{-tL_V^k}\phi^{k,i}](|x|)}{U_k(|x|)}\\
 & =\left[M_{k,i}+o(1)\right]\delta_{0\ell}
 -\left[\frac{N+2A_k}{2}M_{k,i}+o(1)\right]t^{-1}(\partial_r^\ell F_k)(|x|)
 +t^{-2}O(|x|^{4-\ell})\\
 & =\left[M_{k,i}+o(1)\right]\delta_{0\ell}+O(t^{-1}|x|^{2-\ell})\quad\mbox{as}\quad t\to\infty,
\end{split}
\end{equation}
uniformly for $x\in{\bf R}^N$ with $|x|\le\epsilon(1+t)^{\frac{1}{2}}$, 
where $\ell=0,1,2$ and $\delta_{0\ell}$ is the Kronecker symbol. 
Here 
\begin{equation}
\label{eq:2.7}
\begin{split}
  & A_0=A,\quad A_k:=A^+(\lambda_2+\omega_k)\quad\mbox{if}\quad k\ge 1,
 \qquad
 c_{A_k}=\biggr[2^{N+2A_k-1}\Gamma\left(\frac{N+2A_k}{2}\right)\biggr]^{-1/2},\\
 & M_{k,i}:=\frac{c_{A_k}^2}{c_k^2}\int_0^\infty\phi^{k,i}(r)U_k(r)r^{N-1}\,dr
 =\frac{c_{A_k}^2}{c_k^2}\int_{{\bf R}^N}U_k(|y|)Q_{k,i}\left(\frac{y}{|y|}\right)\varphi(y)\,dy,\\
 & F_k(r):=\int_0^r s^{1-N}[U_k(s)]^{-2}\left(\int_0^s \tau^{N-1}U_k(\tau)^2\,d\tau\right)\,ds.
\end{split}
\end{equation}
Here we used 
\begin{equation*}
\begin{split}
\int_{{\bf R}^N}U_k(|y|)Q_{k,i}\left(\frac{y}{|y|}\right)\varphi(y)\,dy
 & =\int_{{\bf S}^{N-1}}Q_{k,i}(\theta)^2\,d\theta\,\int_0^\infty U_k(r)\phi^{k,i}(r)r^{N-1}\,dr\\
 & =\int_0^\infty U_k(r)\phi^{k,i}(r)r^{N-1}\,dr,
\end{split}
\end{equation*}
which follows from the orthonormality of $\{Q_{k,i}\}$ on $L^2({\bf S}^{N-1})$. 
On the other hand, 
it follows from \eqref{eq:1.3}, \eqref{eq:1.4}, \eqref{eq:2.5} and $A_k>A$ that 
\begin{equation}
\label{eq:2.8}
\frac{t^{-\frac{N+2A_k}{2}}U_k(|x|)}{t^{-\frac{N+2A}{2}}U(|x|)}\le Ct^{A-A_k}(1+|x|)^{A_k-A}\le Ct^{-\frac{A_k-A}{2}}=o(1)
\quad\mbox{as}\quad t\to\infty
\end{equation}
uniformly for $x\in{\bf R}^N$ with $|x|\le(1+t)^{\frac{1}{2}}$. 
By Proposition~\ref{Proposition:2.1}, Proposition~\ref{Proposition:2.2}, \eqref{eq:2.6} and \eqref{eq:2.8}, 
for any sufficiently small $\epsilon>0$, 
we have
\begin{equation}
\label{eq:2.9}
[e^{-tL_V^k}\phi^{k,i}](|x|)=o\left([e^{-tL_V}\phi^{0,1}](|x|)\right)\quad\mbox{as}\quad t\to\infty\quad\mbox{if}\quad M_{0,1}>0
\end{equation}
uniformly for $x\in{\bf R}^N$ with $|x|\le\epsilon(1+t)^{\frac{1}{2}}$. 
At the end of this section, we recall the following lemma (see \cite[Lemma~5.1]{IM}).
\begin{lemma}
\label{Lemma:2.1}
Assume condition~{\rm (V)}.
Let $m=1,2,\dots$. Then there exists $C>0$ such that 
\begin{equation}
\label{eq:2.10}
\left|\frac{R_m(x,t)}{U(\min\{|x|,\sqrt{t}\})}\right|\le Ct^{-\frac{N+A^+(\lambda_2+\omega_m)}{2}}U(\sqrt{t})^{-1}
\|\varphi\|_{L^2({\bf R}^N,e^{|x|^2/4}\,dx)}
\end{equation}
for $x\in{\bf R}^N$ and $t>0$. 
\end{lemma}
By Propositions~\ref{Proposition:2.1}, Lemma~\ref{Lemma:2.1}, \eqref{eq:2.6} and \eqref{eq:2.9},  
applying similar arguments as in the proof of \cite[Theorem~1.4]{IM} to \eqref{eq:2.4}, 
we obtain Theorem~\ref{Theorem:1.1}. 
\subsection{Gaussian estimates and the hot spots}
Let $p=p(x,y,t)$ be the fundamental solution generated by $e^{-tL_V}$. 
We first recall the following lemma on upper Gaussian estimates of $p=p(x,y,t)$ 
(see \cite[Theorem~1.3]{IKO}) . 
\begin{proposition}
\label{Proposition:2.3}
Assume condition~{\rm (V)} with $A>-N/2$. 
Then there exists $C>0$ such that 
\begin{equation}
\label{eq:2.11}
0<p(x,y,t)
\le C\,t^{-\frac{N}{2}}
\frac{U(\min\{|x|,\sqrt{t}\})U(\min\{|y|,\sqrt{t}\})}{U(\sqrt{t})^2}
\exp\left(-\frac{|x-y|^2}{Ct}\right)
\end{equation}
for $x$, $y\in{\bf R}^N$ and $t>0$. 
\end{proposition}
By the arguments in Section~2.1 and Proposition~\ref{Proposition:2.3} we have: 
\begin{lemma}
\label{Lemma:2.2}
Let $u$ be a solution of \eqref{eq:1.1} under condition~{\rm(V)}, where $\varphi\in L^2({\bf R}^N,e^{|x|^2/4}\,dx)$. 
Assume that 
$$
\int_{{\bf R}^N}\varphi(y)U(|y|)\,dy>0. 
$$
Then $H(u(t))\not=\emptyset$ for $t>0$. 
Furthermore, there exist $L>0$ and $T>0$ such that 
\begin{equation}
\label{eq:2.12}
H(u(t))\subset B(0,L\sqrt{t})\quad\mbox{for}\quad t\ge T.
\end{equation}
\end{lemma}
{\bf Proof.}
Let $t>0$. 
Since $U$ is a harmonic function for $L_V$, we see that 
$$
U(|x|)=\int_{{\bf R}^N}p(x,y,t)U(|y|)\,dy=\int_{{\bf R}^N}p(y,x,t)U(|y|)\,dy,\qquad x\in{\bf R}^N.
$$
Then the Fubini theorem implies that 
$$
\int_{{\bf R}^N}u(x,t)U(|x|)\,dx=\int_{{\bf R}^N}\int_{{\bf R}^N}p(x,y,t)U(|x|)\varphi(y)\,dy\,dx
=\int_{{\bf R}^N}\varphi(y)U(|y|)\,dy>0. 
$$
Therefore we can find $x_t\in{\bf R}^N$ such that $u(x_t,t)>0$. 
On the other hand, by \eqref{eq:2.11} we can find $R>0$ such that
$$
\sup_{x\in{\bf R}^N\setminus B(0,R)}u(x,t)<u(x_t,t).
$$
This together with \eqref{eq:2.1} implies that $H(u(t))\not=\emptyset$. 

We show \eqref{eq:2.12} in the cases of (S) and (C). 
Since 
$$
|x-y|^2\ge\frac{1}{2}|x|^2-|y|^2,\qquad x,y\in{\bf R}^N,
$$
by \eqref{eq:2.11} we have 
\begin{equation}
\label{eq:2.13}
\begin{split}
 & |u(x,t)|\le Ct^{-\frac{N}{2}}
\left(\int_{B(0,\sqrt{t})}+\int_{{\bf R}^N\setminus B(0,\sqrt{t})}\right)\frac{U(\min\{|y|,\sqrt{t}\})}{U(\sqrt{t})}
e^{-\frac{|x-y|^2}{Ct}}|\varphi(y)|\,dy\\
 & \quad
\le Ct^{-\frac{N+A}{2}}e^{-\frac{|x|^2}{2Ct}}\int_{B(0,\sqrt{t})}e^{\frac{|y|^2}{Ct}}U(|y|)\varphi(y)\,dy
+Ct^{-\frac{N}{2}}e^{-\frac{|x|^2}{2Ct}}\int_{{\bf R}^N\setminus B(0,\sqrt{t})}e^{\frac{|y|^2}{Ct}}\varphi(y)\,dy
\end{split}
\end{equation}
for $x\in{\bf R}^N$ and $t\ge 1$ with $|x|\ge\sqrt{t}$. 
Recalling that $\varphi\in L^2({\bf R}^N,e^{|x|^2/4}\,dx)$, 
by the Cauchy-Schwarz inequality we see that
\begin{equation}
\label{eq:2.14}
\begin{split}
\int_{{\bf R}^N\setminus B(0,\sqrt{t})}e^{\frac{|y|^2}{Ct}}|\varphi(y)|\,dy
 & \le\left(\int_{{\bf R}^N\setminus B(0,\sqrt{t})}e^{\frac{2|y|^2}{Ct}}e^{-\frac{|y|^2}{4}}\,dy\right)^{\frac{1}{2}}
\left(\int_{{\bf R}^N}e^{\frac{|y|^2}{4}}|\varphi(y)|^2\,dy\right)^{\frac{1}{2}}\\
 & \le C\left(\int_{{\bf R}^N\setminus B(0,\sqrt{t})}e^{-\frac{|y|^2}{8}}\,dy\right)^{\frac{1}{2}}
 \le Ce^{-\frac{t}{C}}
\end{split}
\end{equation}
for all sufficiently large $t$. 
Then, 
for any $\epsilon>0$, by \eqref{eq:1.4}, \eqref{eq:2.13} and \eqref{eq:2.14}
we can find constants $L\ge 1$ such that 
\begin{equation}
\label{eq:2.15}
\sup_{x\in{\bf R}^N \backslash B(0,L\sqrt{t})}|u(x,t)|\le \epsilon t^{-\frac{N+A}{2}}
\end{equation}
for all sufficiently large $t$. 
On the other hand, by Theorem~\ref{Theorem:1.1}~(a) we see that 
\begin{equation}
\label{eq:2.16}
 \liminf_{t\to\infty}t^{-\frac{N+A}{2}}\max_{y\in\partial B(0,1)}u(\sqrt{t}y,t)\ge CM(\varphi)>0. 
\end{equation}
Taking a sufficiently small $\epsilon>0$ if necessary, 
we deduce from \eqref{eq:2.15} and \eqref{eq:2.16} that $H(u(t))\subset B(0,L\sqrt{t})$ for all sufficiently large $t$. 
Thus \eqref{eq:2.12} holds in the cases of (S) and (C). 
Similarly, we can prove \eqref{eq:2.12} in the case of ($\mbox{S}_*$).  
Thus Lemma~\ref{Lemma:2.2} follows. 
$\Box$
\vspace{5pt}

By Theorem~\ref{Theorem:1.1} and Lemma~\ref{Lemma:2.2} we have:
\begin{theorem}
\label{Theorem:2.1}
Let $L_V$ be a nonnegative Schr\"odinger operator under condition~{\rm (V)} with $\lambda_1<0$.  
Let $u$ be a solution of \eqref{eq:1.1} such that 
\begin{equation}
\label{eq:2.17}
\varphi\in L^2({\bf R}^N,e^{|x|^2/4}\,dx),
\qquad
M(\varphi):=\int_{{\bf R}^N}\varphi(y)U(|y|)\,dy>0. 
\end{equation}
Then $H(u(t))=\{0\}$ for all sufficiently large $t>0$. 
\end{theorem}
In the case of $\lambda_1<0$, Theorem~\ref{Theorem:1.1} together with
\eqref{eq:1.3} implies that $A^+(\lambda_1)<0$, $U(r)\thicksim r^{A^+(\lambda_1)}$ as $r\to 0$ 
and $u(0,t)~{=+\infty}$ for all sufficiently large $t$. 
Thus Theorem~\ref{Theorem:2.1} follows. 
\section{Large time behavior of the hot spots}
We study the large time behavior of the hot spots for problem~\eqref{eq:1.1} in the case of $\lambda_1\ge 0$. 
Set 
$$
\Pi:=\biggr\{r\in[0,\infty)\,:\,U(r)=\sup_{\tau\in[0,\infty)}U(\tau)\biggr\},
\qquad
\Xi(\varphi):=\int_{{\bf R}^N}\varphi(y)U_1(|y|)\frac{y}{|y|}dy.
$$
Throughout this section we use the same notation as in Section~2. For the reader's convenience, 
we give a correspondence table of our theorems. 
We recall $A=A^+(\lambda_2)$ in the cases of (S) and (S$_{*}$) and $A=A^-(\lambda_2)$ in the case of (C). 
\begin{itemize}
  \item[(I)\,] $\lambda_1<0$ (see Theorem~\ref{Theorem:2.1});
  \item[(II)] $\lambda_1\ge 0$ 
  \begin{itemize}
  \item[(1)] $A>0$ (see Theorem~\ref{Theorem:3.1});
  \item[(2)] $A=0$
  	\begin{itemize}
	\item[(a)] $\Pi=\emptyset$ and $N=2$ (see Theorems~\ref{Theorem:3.2} and \ref{Theorem:3.3}); 
	\item[(b)] $\Pi=\emptyset$ and $N\ge 3$ (see Theorem~\ref{Theorem:3.3} 
	and Corollaries~\ref{Corollary:3.1} and \ref{Corollary:3.2});
	\item[(c)] $\Pi\not=\emptyset$ (see Theorems~\ref{Theorem:3.6} and \ref{Theorem:3.7});
  	\end{itemize}
  \item[(3)] $A<0$ (see Theorems~\ref{Theorem:3.4} and \ref{Theorem:3.5}). 
  \end{itemize}
\end{itemize}
We consider the case of $A>0$. 
Theorem~\ref{Theorem:3.1} is proved by the same arguments as in \cite[Section~4]{IK02} 
with the aid of the results in Section~2. See also Theorem~\ref{Theorem:1.1} and \cite[Theorem~1.2]{IK02}. 
\begin{theorem}
\label{Theorem:3.1}
Let $L_V$ be a nonnegative Schr\"odinger operator under condition~{\rm (V)} with $\lambda_1\ge 0$ and $A>0$.  
Assume \eqref{eq:2.17} and let $u$ be a solution of \eqref{eq:1.1}. 
Then the following holds:
\begin{itemize}
  \item[{\rm (a)}] 
  $\displaystyle{\lim_{t\to\infty}\sup_{x\in H(u(t))}\left|\,t^{-\frac{1}{2}}|x|-\sqrt{2A}\,\right|=0}$;
  \item[{\rm (b)}] 
  Assume that $\Xi(\varphi)\not=0$. 
  Then there exist a constant $T>0$ and a curve $x=x(t)\in C^1([T,\infty):{\bf R}^N)$ such that 
  $H(u(t))=\{x(t)\}$ for $t\ge T$ and 
  $$
  \lim_{t\to\infty}\frac{x(t)}{|x(t)|}=\frac{\Xi(\varphi)}{|\Xi(\varphi)|}. 
  $$
  \end{itemize}
\end{theorem}
Secondly, we consider the case where $A=0$ and $\Pi=\emptyset$.  
Theorems~\ref{Theorem:3.2} and \ref{Theorem:3.3} 
are obtained by the same arguments as in 
\cite[Section~4]{IK04} and \cite[Sections~4 and 5]{IK03}, respectively, with the aid of the results in Section~2. 
See also \cite[Theorem~1.2]{IK04} and \cite[Theorem~1.2]{IK03}. 
We remark that, in the case of $N=2$, $A=0$ if and only if $\lambda_2=0$. 
Furthermore, if $L_V$ is subcritical, then $U(r)\asymp \log r$ as $r\to\infty$ and $\Pi=\emptyset$. 
\begin{theorem}
\label{Theorem:3.2} 
Let $N=2$ and $L_V$ be a subcritical Schr\"odinger operator under condition~{\rm (V)} 
with $\lambda_1\ge 0$ and $A=0$. 
Assume \eqref{eq:2.17} and let $u$ be a solution of \eqref{eq:1.1}. 
Then
$$
\lim_{t\to\infty}\sup_{x\in H(u(t))}\left|t^{-1}(\log t)|x|-2\right|=0. 
$$
Furthermore, assertion~{\rm (b)} of Theorem~{\rm\ref{Theorem:3.1}} holds. 
\end{theorem}
\begin{theorem}
\label{Theorem:3.3} 
Let $L_V$ be a nonnegative Schr\"odinger operator under condition~{\rm (V)} with 
$\lambda_1\ge 0$, $A=0$ and $\Pi=\emptyset$. 
Assume that $L_V$ is critical if $N=2$. 
Assume \eqref{eq:2.17} and let $u$ be a solution of \eqref{eq:1.1}. 
Then 
$$
\lim_{t\to\infty}\sup_{x\in H(u(t))}\left|\frac{tU'(|x|)}{c_*|x|}-\frac{1}{2}\right|=0.
$$
Furthermore, assertion~{\rm (b)} of Theorem~{\rm\ref{Theorem:3.1}} holds. 
\end{theorem}
Let $\lambda_1\ge 0$ and 
$$
\Gamma_k(r):=r^k\int_0^rs^{1-N-2k}\left(\int_0^s\tau^{N+k-1}V(\tau)U_k(\tau)\,d\tau\right)\,ds,
\qquad k=0,1,2,\dots.
$$
Since $\tilde{U}_k:=U_k-\Gamma_k$ satisfies 
$$
\tilde{U}_k''+\frac{N-1}{r}\tilde{U}_k'-\frac{\omega_k}{r^2}\tilde{U}_k=0\quad\mbox{in}\quad(0,\infty),
\qquad
\tilde{U}_k(r)\asymp r^{A^+(\lambda_1+\omega_k)}\quad\mbox{as}\quad r\to 0, 
$$
by the uniqueness of the solution of \eqref{eq:1.3} with $V$ replaced by $\omega_k r^{-2}$ 
we see that 
\begin{equation}
\label{eq:3.1}
U_k(r)=r^k+\Gamma_k(r)\quad\mbox{if}\quad\lambda_1=0,\qquad
U_k(r)=\Gamma_k(r)\quad\mbox{if}\quad\lambda_1>0.
\end{equation}
In particular, in the case of $A=0$, we have 
$$
c_*\equiv \lim_{r\to\infty}U(r)=
\left\{
\begin{array}{ll}
\Gamma_0(\infty)+1 & \mbox{if}\quad\lambda_1=0,\vspace{3pt}\\
\Gamma_0(\infty) & \mbox{if}\quad\lambda_1>0. 
\end{array}
\right.
$$
As a corollary of Theorem~\ref{Theorem:3.3}, we have the following result, 
which revises \cite[Corollary~1.1]{IK03} and \cite[Remark~1.1]{IK04}.
\begin{corollary}
\label{Corollary:3.1}
Assume the same conditions as in Theorem~{\rm\ref{Theorem:3.3}} with $\lambda_1=0$. 
Furthermore, assume that $V(r)\thicksim \mu r^{-d}$ as $r\to\infty$ for some $\mu\not=0$ and $d>2$.
\begin{itemize}
  \item[{\rm(a)}] Let $\mu>0$. Then 
  $$
  |x|=
  \left\{
  \begin{array}{ll}
  \displaystyle{\left(\frac{2\mu t}{(\Gamma_0(\infty)+1)(N-d)}\right)^{\frac{1}{d}}(1+o(1))}
   & \mbox{if}\quad 2<d<N,\vspace{5pt}\\
   \displaystyle{\left(\frac{2\mu\,t\log t}{(\Gamma_0(\infty)+1)N}\right)^{\frac{1}{N}}(1+o(1))}
   & \mbox{if}\quad d=N,\vspace{5pt}\\
    \displaystyle{\left(\frac{2\Lambda t}{\Gamma_0(\infty)+1}\right)^{\frac{1}{N}}(1+o(1))}
   & \mbox{if}\quad d>N,\quad \Lambda>0,
  \end{array}
  \right.
  $$
  as $t\to\infty$ uniformly for $x\in H(u(t))$. Here 
  $\displaystyle{\Lambda:=\int_0^\infty \tau^{N-1}V(\tau)U(\tau)\,d\tau}$. 
 \item[{\rm (b)}] 
 Let $\mu<0$ and $d>N$. Then 
 $$
 |x|=
 \left\{
 \begin{array}{ll}
 \displaystyle{\left(\frac{2\Lambda t}{\Gamma_0(\infty)+1}\right)^{\frac{1}{N}}(1+o(1))} & \mbox{if}\quad\Lambda>0,\vspace{5pt}\\
 \displaystyle{\left(\frac{2|\mu|t}{(\Gamma_0(\infty)+1)(d-N)}\right)^{\frac{1}{d}}(1+o(1))} & \mbox{if}\quad\Lambda=0,
 \end{array}
 \right.
 $$
 as $t\to\infty$ uniformly for $x\in H(u(t))$.  
\end{itemize}
\end{corollary}
Furthermore, we have: 
\begin{corollary}
\label{Corollary:3.2}
Assume the same conditions as in Theorem~{\rm\ref{Theorem:3.3}} with $\lambda_1>0$. 
Then the same assertions of Corollary~{\rm\ref{Corollary:3.1}} holds with $\Gamma_0(\infty)+1$ replaced by $\Gamma_0(\infty)$. 
\end{corollary}
\begin{remark}
\label{Remark:3.1}
Assume the same conditions as in Theorem~{\rm\ref{Theorem:3.3}}. 
Let $V(r)\thicksim \mu r^{-d}$ as $r\to\infty$ for some $\mu\not=0$ and $d>2$. 
\newline
{\rm (i)} Consider the case where $\mu>0$ and  $d>N$. 
Since $U(r)\thicksim c_*>0$ as $r\to\infty$, $\Lambda$ can be defined. 
If $\Lambda\le 0$ and $\mu>0$, 
then it follows from \eqref{eq:3.1} that
$$
U'(r)=r^{1-N}\int_0^r \tau^{N-1}V(\tau)U(\tau)\,d\tau
=r^{1-N}\left[\Lambda-\int_r^\infty \tau^{N-1}V(\tau)U(\tau)\,d\tau\right]<0
$$
for all sufficiently large $r>0$. This implies that $\Pi\not=\emptyset$. 
\newline
{\rm (ii)} 
Consider the case where $\mu<0$. 
By \eqref{eq:3.1} we see that $r^{N-1}U'(r)\to -\infty$ as $r\to\infty$ if $2<d\le N$.  
Similarly, if $d>N$ and $\Lambda<0$, 
then $U'(r)<0$ for all sufficiently large $r>0$. 
In the both cases, it follows that $\Pi\not=\emptyset$. 
\end{remark}
\vspace{3pt}

Next we study the large time behavior of the hot spots in the case where $\lambda_1\ge 0$ and $A<0$. 
It follows from $A<0$ that  $U(r)\to 0$ as $r\to\infty$ and $\Pi\not=\emptyset$. 
\begin{theorem}
\label{Theorem:3.4}
Let $L_V$ be a nonnegative Schr\"odinger operator under condition~{\rm (V)} with 
$\lambda_1\ge 0$ and $A<0$. 
Assume \eqref{eq:2.17} and let $u$ be a solution of \eqref{eq:1.1}. 
Then  
\begin{equation}
\label{eq:3.2}
\lim_{t\to\infty}\sup_{x\in H(u(t))}\left|\,|x|-\min\Pi\,\right|=0. 
\end{equation}
Furthermore, if $\Xi(\varphi)\not=0$, then 
\begin{equation}
\label{eq:3.3}
\lim_{t\to\infty}\sup_{x\in H(u(t))}\left|\,x-\min\Pi\,\frac{\Xi(\varphi)}{|\Xi(\varphi)|}\,\right|=0.
\end{equation}
\end{theorem}
{\bf Proof.} 
For any $\epsilon>0$, 
by Theorem~\ref{Theorem:1.1} with $A<0$ and Lemma~\ref{Lemma:2.2} we see that 
\begin{equation}
\label{eq:3.4}
H(u(t))\subset B(0,\epsilon\sqrt{t})
\end{equation}
for all sufficiently large $t$. 

We consider the cases of (S) and (C). 
In the case of $\Xi(\varphi)\not=0$ 
we can assume, without loss of generality, that $\Xi(\varphi)=(|\Xi(\varphi)|,0,\dots,0)$. 
By \eqref{eq:2.7} we have
\begin{equation}
\label{eq:3.5}
M_{0,1}>0,\qquad
M_{1,i}=\frac{c_{A_1}^2}{c_1^2}\Xi_i(\varphi)=\frac{c_{A_1}^2}{c_1^2}|\Xi(\varphi)|\delta_{1,i},
\quad i=1,\dots,N,
\end{equation}
where $\delta_{1,i}$ is the Kronecker symbol.
Let $\epsilon>0$ be sufficiently small. 
By Lemma~\ref{Lemma:2.1}, \eqref{eq:1.4}, \eqref{eq:2.3}, \eqref{eq:2.6} and \eqref{eq:2.9}
we take a sufficiently large $m\in\{1,2,\dots\}$ so that 
\begin{equation}
\label{eq:3.6}
\begin{split}
 & t^{\frac{N+2A}{2}}\frac{u_0(x,t)}{q_*U(|x|)}=[M_{0,1}+o(1)]
-\left[\frac{N+2A}{2}M_{0,1}+o(1)\right]t^{-1}F(|x|)\\
 & \qquad\qquad\qquad\qquad\qquad\qquad\qquad\qquad\qquad\qquad\qquad\qquad\qquad\quad
 +t^{-2}O(|x|^4)\\
 & \qquad\qquad\qquad
 =[M_{0,1}+o(1)]+t^{-1}O(|x|^2),\\
 & u_{1,i}(x,t)=q_N[M_{1,i}+o(1)]t^{-\frac{N+2A_1}{2}}U_1(|x|)\frac{x_i}{|x|}
+t^{-\frac{N+2A_1}{2}}U_1(|x|)O(t^{-1}|x|^2)\\
 & \qquad\quad\,\,\,\,
 =O(t^{-\frac{N+2A_1}{2}}(1+|x|)^{A_1}),\\
 & R_2(x,t)=\sum_{k=2}^{m-1}\sum_{i=1}^{\ell_k}u_{k,i}(x,t)+R_m(x,t)
=O\left(t^{-\frac{N+2A_2}{2}}(1+|x|)^{A_2}\right),
\end{split}
\end{equation}
as $t\to\infty$ uniformly for $x\in{\bf R}^N$ with $|x|\le\epsilon(1+t)^{\frac{1}{2}}$, where $i=1,\dots,N$. 
Let $\nu$ be a sufficiently small positive constant. 
Since $A<0$, 
we can find $R>0$ such that 
\begin{equation}
\label{eq:3.7}
\begin{split}
 & \left|u_0(x,t)\right|\le Ct^{-\frac{N+2A}{2}}U(|x|)\le \nu t^{-\frac{N+2A}{2}},\\
 & \left|u_{1,i}(x,t)\right|\le \nu t^{-\frac{N+2A}{2}},
 \qquad |R_2(x,t)|\le \nu t^{-\frac{N+2A}{2}},
\end{split}
\end{equation}
for $x\in{\bf R}^N$ and all sufficiently large $t>0$ with $|x|\le\epsilon(1+t)^{\frac{1}{2}}$, where $i=1,\dots,N$. 
On the other hand, Theorem~\ref{Theorem:1.1} implies that 
\begin{equation}
\label{eq:3.8}
\liminf_{t\to\infty}t^{\frac{N+2A}{2}}\sup_{x\in{\bf R}^N}u(x,t)>0. 
\end{equation}
By \eqref{eq:3.4}, \eqref{eq:3.7} and \eqref{eq:3.8} 
we can find $R>0$ such that 
\begin{equation}
\label{eq:3.9}
H(u(t))\subset B(0,R)
\end{equation}
for all sufficiently large $t$. 

It follows from \eqref{eq:1.2} and $\lambda_2<0$ that 
\begin{equation}
\label{eq:3.10}
A+1\le A^+(\lambda_2)+1<A^+(\lambda_2+\omega_1)=A_1.
\end{equation}
By \eqref{eq:2.4}, \eqref{eq:3.6} and \eqref{eq:3.10} we have
\begin{equation}
\label{eq:3.11}
\begin{split}
 & t^{\frac{N+2A}{2}}u(x,t)\\
 & =q_*[M_{0,1}+o(1)]U(|x|)-\left[\frac{N+2A}{2}M_{0,1}+o(1)\right]t^{-1}U(|x|)F(|x|)+o(t^{-1})\\
 & =q_*[M_{0,1}+o(1)]U(|x|)+O(t^{-1})
\end{split}
\end{equation}
as $t\to\infty$ uniformly for $x\in B(0,R)$. 
Since $F$ is strictly monotone increasing in $(0,\infty)$, 
by \eqref{eq:3.5}, \eqref{eq:3.9} and \eqref{eq:3.11} we obtain \eqref{eq:3.2} and \eqref{eq:3.3}. 
Therefore Theorem~\ref{Theorem:3.4} follows in the cases of (S) and (C). 
Similarly, Theorem~\ref{Theorem:3.4} also follows in the case of $(\mbox{S}_*)$. 
Thus the proof is complete.
$\Box$\vspace{3pt}
\newline
We give sufficient conditions for the hot spots to 
consist of only one point and to move along a smooth curve. 
We denote by $\nabla^2 f$ the Hessian matrix of a function $f$. 
For any real symmetric $N\times N$ matrix $M$, 
by $M\ge 0$ and $M\le 0$ we mean that $M$ is positive semi-definite and negative semi-definite, respectively. 
\begin{theorem}
\label{Theorem:3.5}
Let $L_V$ be a nonnegative Schr\"odinger operator under condition~{\rm (V)} with 
$\lambda_1\ge 0$ and $A<0$. 
Assume \eqref{eq:2.17} and let $u$ be a solution of \eqref{eq:1.1}. 
Let $x_*\in{\bf R}^N$ be such that $|x_*|\in\Pi$ and 
$$
\lim_{t\to\infty}\sup_{x\in H(u(t))}|x-x_*|=0.
$$
Then there exist a constant $T>0$ and a curve $x=x(t)\in C^1([T,\infty):{\bf R}^N)$ such that 
$H(u(t))=\{x(t)\}$ for $t\ge T$ in the following cases: 
\begin{itemize}
  \item[{\rm (a)}] 
  $|x_*|=0$, $V\in C^\gamma([0,\infty))$ for some $\gamma\in(0,1)$ and 
  $\nabla^2 U(|x|)\le 0$ in a neighborhood of $x=0$;
  \item[{\rm (b)}] 
  $|x_*|>0$, $U''\le 0$ in a neighborhood of $r=|x_*|$ and $\Xi(\varphi)\not=0$. 
\end{itemize}
\end{theorem}
{\bf Proof.}
We consider the cases of (S) and (C). 
Let $r_*:=|x_*|$ and $\epsilon>0$. 
The proof is divided into the following four cases:
\begin{equation*}
\begin{array}{llll}
{\rm (I)}& r_*=0,\quad U''(0)<0;
\qquad
 & {\rm (II)} & r_*=0,\quad U''(0)=0;\vspace{3pt}\\
{\rm (III)} & r_*>0,\quad U''(r_*)<0;
\qquad
 & {\rm (IV)} & r_*>0,\quad U''(r_*)=0.
\end{array}
\end{equation*}
We consider case~(I). 
Since $U'(0)=0$ and $U''(0)<0$, by \eqref{eq:2.3} and \eqref{eq:2.6} we can find $\eta_1>0$ such that
\begin{equation}
\label{eq:3.12}
-q_*^{-1}t^{\frac{N+2A}{2}}(\nabla^2u_0)(x,t)
=-[M_{0,1}+o(1)](\nabla^2 U)(|x|)+O(t^{-1})
\ge -\frac{M_{0,1}U''(0)}{2}I_N-\epsilon I_N
\end{equation}
for $x\in B(0,\eta_1)$ and all sufficiently large $t$, where $I_N$ is the $N$-dimensional identity matrix. 
On the other hand, 
by condition~(a), \eqref{eq:2.6} and \eqref{eq:2.10} 
we apply the parabolic regularity theorems to see that 
$u_{1,i}$, $R_2\in C^{2,\gamma;1,\gamma/2}({\bf R}^N\times(0,\infty))$ and 
\begin{equation}
\label{eq:3.13}
\|\nabla^2 u_{1,i}\|_{L^\infty(B(0,\eta_1))}+\|\nabla^2 R_2\|_{L^\infty(B(0,\eta_1))}=O(t^{-\frac{N+2A_1}{2}})
\end{equation}
for all sufficiently large $t$, where $i=1,\dots,N$.
Since $\epsilon$ is arbitrary, by \eqref{eq:3.12} and \eqref{eq:3.13} we see that 
$-(\nabla^2 u)(x,t)$ is positive definite in $B(0,\eta_1)$ for all sufficiently large $t>0$. 
Then Theorem~\ref{Theorem:3.5} in case~(I) follows from the implicit function theorem. 

Consider case (II). 
By condition~(a), \eqref{eq:2.3} and \eqref{eq:2.6} we have 
\begin{equation}
\label{eq:3.14}
\begin{split}
 & -q_*^{-1}t^{\frac{N+2A}{2}+1}(\nabla^2 u_0)(x,t)\\
 & =-t[M_{0,1}+o(1)]\nabla^2 U(|x|)
+\left[\frac{N+2A}{2}M_{0,1}+o(1)\right]\nabla^2[UF_0](|x|)+O(t^{-1})\\
 & \ge\left[\frac{N+2A}{2}M_{0,1}+o(1)\right]\nabla^2[UF_0](|x|)+O(t^{-1})
\end{split}
\end{equation}
in a neighborhood of $x=0$ and all sufficiently large $t$. 
On the other hand, it follows from \eqref{eq:1.3} and \eqref{eq:2.7} that 
\begin{equation}
\label{eq:3.15}
\begin{split}
F_0''(r) & =\left\{(1-N)r^{-N}[U(r)]^{-2}-2r^{1-N}[U(r)]^{-3}U'(r)\right\}\int_0^r \tau^{N-1}U(\tau)^2\,d\tau+1\\
 & \to\frac{1-N}{N}+1=\frac{1}{N}\quad\mbox{as}\quad r\to 0.
\end{split}
\end{equation}
This implies that $[UF_0]''(0)=1/N$. 
Therefore, by \eqref{eq:3.14} and \eqref{eq:3.15} we can find $\eta_2>0$ such that  
\begin{equation}
\label{eq:3.16}
-q_*^{-1}t^{\frac{N+2A}{2}+1}(\nabla^2 u_0)(x,t)
\ge\frac{N+2A}{4N}M_{0,1}I_N-\epsilon I_N
\end{equation}
for $x\in B(0,\eta_2)$ and all sufficiently large $t$. 
Similarly to case~(I), since $\epsilon$ is arbitrary, 
by \eqref{eq:2.4}, \eqref{eq:3.13} and \eqref{eq:3.16} we see that 
$-(\nabla^2 u)(x,t)$ is positive definite in $B(0,\eta_2)$ for all sufficiently large $t>0$. 
Similarly to case~(I), 
Theorem~\ref{Theorem:3.5} in case~(II) follows from the implicit function theorem. 

Consider case~(III). 
By Theorem~\ref{Theorem:3.4} 
we can assume, without loss of generality, that $x_*=(r_*,0,\dots,0)$. Then $M_{1,1}>0$ and $M_{1,i}=0$ for $i\in\{2,\dots,N\}$. 
Let $\theta_\alpha:=x_\alpha/|x|$ for $\alpha=1,\dots,N$. 
Then $(r,\theta_2,\dots,\theta_N)$ gives a local coordinate of ${\bf R}^N$ in a neighborhood of $x_*$. 
We study the large time behavior of $\tilde{\nabla}^2u$ in a neighborhood of $x_*$, 
where $\tilde{\nabla}:=(\partial_r,\partial_{\theta_2},\dots,\partial_{\theta_N})$. 
Since $U''(r_*)<0$, 
similarly to \eqref{eq:3.12}, we can find $\eta_3>0$ such that 
\begin{equation}
\label{eq:3.17}
-q_*^{-1}t^{\frac{N+2A}{2}}(\partial_r^2 u_0)(x,t)
=-[M_{0,1}+o(1)](\partial_r^2 U)(|x|)+O(t^{-1})
\ge -\frac{M_{0,1}U''(r_*)}{2}
\end{equation}
for $x\in B(x_*,\eta_3)$ and all sufficiently large $t$. 
Furthermore, 
\begin{equation}
\label{eq:3.18}
(\partial_r\partial_{\theta_\alpha}u_0)(x,t)=(\partial_{\theta_\alpha}\partial_{\theta_\beta}u_0)(x,t)=0
\end{equation}
for $x\in B(x_*,\eta_3)$ and all sufficiently large $t$, where $\alpha$, $\beta\in\{2,\dots,N\}$. 

On the other hand, 
by \eqref{eq:2.3}, \eqref{eq:2.6} and \eqref{eq:3.6} we have 
\begin{equation}
\label{eq:3.19}
q_N^{-1}(\tilde{\nabla}^2 u_{1,i})(x,t)
=[M_{1,i}+o(1)]t^{-\frac{N+2A_1}{2}}\tilde{\nabla}^2[U_1(|x|)\theta_i]
+O(t^{-\frac{N+2A_1}{2}-1})
\end{equation}
for $x\in B(x_*,\eta_3)$ and all sufficiently large $t$. 
Since 
\begin{equation}
\label{eq:3.20}
\theta_1=\biggr(1-\sum_{\alpha=2}^N\theta_\alpha^2\biggr)^{1/2},\quad
\frac{\partial\theta_1}{\partial\theta_\alpha}=-\theta_1^{-1}\theta_\alpha,\quad
\frac{\partial^2\theta_1}{\partial\theta_\alpha\partial\theta_\beta}
=-\delta_{\alpha\beta}\theta_1^{-1}
-\theta_1^{-3}\theta_\alpha\theta_\beta,
\end{equation}
for $\alpha$, $\beta\in\{2,\dots,N\}$, 
combining $M_{1,1}=|\Xi(\varphi)|>0$ and $U'(r_*)=0$, 
we can find $\eta_4>0$ and $C>0$ such that  
\begin{equation}
\label{eq:3.21}
\begin{split}
 & -t^{\frac{N+2A_1}{2}}(\partial_r^2u_{1,1})(x,t)\ge -C,\\
 & -t^{\frac{N+2A_1}{2}}(\partial_{\theta_\alpha}\partial_{\theta_\beta}u_{1,1})(x,t)
\ge\frac{q_NM_{1,1}}{2}U_1(r_*)\delta_{\alpha\beta}-\epsilon,\\
 & -t^{\frac{N+2A_1}{2}}(\partial_r\partial_{\theta_j}u_{1,1})(x,t)
\ge- C|\theta_\alpha|U_1'(r)+O(t^{-1})\ge-\epsilon,
\end{split}
\end{equation}
for $x\in B(x_*,\eta_4)$ and all sufficiently large $t$. 
Furthermore, 
for $i=2,\dots,N$, 
it follows that $M_{1,i}=0$ and we have 
\begin{equation}
\label{eq:3.22}
(\tilde{\nabla}^2u_{1,i})(x,t)=o(t^{-\frac{N+2A_1}{2}})
\end{equation}
for $x\in B(x_*,\eta_4)$ and all sufficiently large $t$. 
Similarly to \eqref{eq:3.13}, 
by \eqref{eq:3.6} we apply the parabolic regularity theorems to obtain 
\begin{equation}
\label{eq:3.23}
(\tilde{\nabla}^2R_2)(x,t)=O(t^{-\frac{N+2A_2}{2}})
\end{equation}
for $x\in B(x_*,\eta_4)$ and all sufficiently large $t$. 
On the other hand, $A_1>A+1$ holds by $A<0$. 
Then, by \eqref{eq:3.17}, \eqref{eq:3.18}, \eqref{eq:3.19}, \eqref{eq:3.21}, \eqref{eq:3.22} and \eqref{eq:3.23}
we see that 
$-(\tilde{\nabla}^2u)(x,t)$ is positive definite in a neighborhood of $x_*=(r_*,0)$ for all sufficiently large $t>0$. 
Therefore Theorem~\ref{Theorem:3.5} in case~(III) follows from the implicit function theorem. 

It remains to consider case~(IV). 
Similarly to case~(III), without loss of generality, we can assume that $\Xi(\varphi)/|\Xi(\varphi)|=(1,0,\dots,0)$. 
It follows from $U'(r_*)=U''(r_*)=0$ and $r_*\in\Pi$ that 
\begin{equation*}
\begin{split}
[UF_0]''(r_*) & =U(r_*)F_0''(r_*)\\
 & =U(r_*)\left[(1-N)r^{-N}[U(r_*)]^{-2}\int_0^{r_*} \tau^{N-1}U(\tau)^2\,d\tau+1\right]\\
 & \ge U(r_*)\left[(1-N)r_*^{-N}[U(r_*)]^{-2}\int_0^{r_*}\tau^{N-1}U(r_*)^2\,d\tau+1\right]
 =\frac{1}{N}U(r_*)>0. 
\end{split}
\end{equation*}
Then, by condition~(b) we have 
\begin{equation}
\label{eq:3.24}
\begin{split}
 & -q_*^{-1}t^{\frac{N+2A}{2}+1}(\partial_r^2 u_0)(x,t)\\
 & =-[M_{0,1}+o(1)]t(\partial_r^2 U)(|x|)+\left[\frac{N+2A}{2}M_{0,1}+o(1)\right]\partial_r^2[UF](|x|)+O(t^{-1})\\
 & \ge \left[\frac{N+2A}{2}M_{0,1}+o(1)\right]\partial_r^2[UF](|x|)+O(t^{-1})
\ge\frac{N+2A}{4N}M_{0,1}U(r_*)>0
\end{split}
\end{equation}
in a neighborhood of $x_*=(r_*,0,\dots,0)$. 
Furthermore, by the same argument as in case~(III) 
we obtain \eqref{eq:3.18}, \eqref{eq:3.21}, \eqref{eq:3.22} and \eqref{eq:3.23}. 
Therefore, since $A_1> A+1$, 
we see that $-(\tilde{\nabla}^2u)(x,t)$ is positive definite in a neighborhood of $x_*=(r_*,0)$ for all sufficiently large $t>0$. 
Therefore Theorem~\ref{Theorem:3.5} in case~(III) follows from the implicit function theorem. 
Thus Theorem~\ref{Theorem:3.5} follows in the cases of (S) and (C). 
Similarly, Theorem~\ref{Theorem:3.5} also follows in case ({\rm $S_*$). 
Therefore the proof of Theorem~\ref{Theorem:3.5} is complete. 
$\Box$
\vspace{5pt}

Finally we study the large time behavior of the hot spots in the cases where $\lambda_1\le 0$, $A=0$ 
and $\Pi\not=\emptyset$. 
\begin{theorem}
\label{Theorem:3.6}
Let $L_V$ be a nonnegative Schr\"odinger operator under condition~{\rm (V)} with 
$\lambda_1\ge 0$,  $A=0$ and $\Pi\not=\emptyset$. 
Assume \eqref{eq:2.17} and let $u$ be a solution of \eqref{eq:1.1}. 
Then there exists $R>0$ such that 
\begin{equation}
\label{eq:3.25}
H(u(t))\subset B(0,R)
\end{equation}
for all sufficiently large $t$. 

Let $x_*$ be an accumulating point of $H(u(t))$ as $t\to\infty$. 
If $\Xi(\varphi)\not=0$, then $x_*=|x_*|\Xi(\varphi)$. Furthermore, $r_*:=|x_*|$ is a maximum point of 
$$
S(r):=-\frac{N}{c_*^2}M(\varphi)U(r)F_0(r)
+\frac{1}{c_1^2}|\Xi(\varphi)|\,U_1(r)\quad\mbox{on}\quad\Pi.
$$
\end{theorem}
\begin{remark}
\label{Remark:3.2}
Assume the same conditions as in Theorem~{\rm\ref{Theorem:3.6}}. 
It follows from $\lambda_2=0$ that $A=0$ and $A_1=1$. 
By \eqref{eq:2.2} and \eqref{eq:2.5} we see that 
\begin{equation}
\label{eq:3.26}
F_0(r)\thicksim(2N)^{-1}r^2,\qquad U(r)\thicksim c_*,\qquad U_1(r)\thicksim c_1r,
\end{equation}
as $r\to\infty$. Then $S(x)\to-\infty$ as $|x|\to\infty$ and the maximum point of $S$ on $\Pi$ exists. 
\end{remark}
{\bf Proof.} 
Let $\epsilon$ be a sufficiently small positive constant. 
Similarly to \eqref{eq:3.4}, 
by Theorem~\ref{Theorem:1.1} with $A=0$ we see that 
\begin{equation}
\label{eq:3.27}
H(u(t))\subset B(0,\epsilon\sqrt{t})
\end{equation}
for all sufficiently large $t$. 
On the other hand, similarly to the proof of Theorem~\ref{Theorem:3.4}, 
we have 
\begin{equation}
\label{eq:3.28}
\begin{split}
t^{\frac{N}{2}} u(x, t)
& =q_*[M_{0, 1} + o(1)]U(|x|) -q_*\left[ \frac{N}{2} M_{0, 1} +o(1) \right] t^{-1} (UF_0)(|x|)\\
&\qquad  + t^{-2}O(|x|^4U(|x|))+ \sum^N_{i=1}  [M_{1, i}+o(1)] t^{-1}  U_1(|x|) Q_{1, i} \left( \frac{x}{|x|} \right) \\
& \qquad\qquad\qquad\qquad\qquad\qquad\qquad\qquad+t^{-2} U_1(|x|) O(|x|^2) +o(t^{-1})
\end{split}
\end{equation}
as $t\to\infty$ uniformly for $x\in{\bf R}^N$ with $|x|\le\epsilon t^{1/2}$. 
Since $A_1=A^+(\omega_1)=1$, it follows from \eqref{eq:2.3}, \eqref{eq:2.7} and \eqref{eq:2.17} that 
\begin{equation*}
\begin{split}
q_*M_{0,1} & =q_*^2\frac{1}{c_*^2}\biggr[2^{N-1}\Gamma\biggr(\frac{N}{2}\biggr)\biggr]^{-1}M(\varphi),\\
\sum_{i=1}^NM_{1, i}Q_{1, i} \left( \frac{x}{|x|} \right)
 & =\frac{1}{c_1^2}\biggr[2^{N+1}\Gamma\biggr(\frac{N+2}{2}\biggr)\biggr]^{-1}q_N^2
 \biggr(\frac{x}{|x|}\cdot\int_{{\bf R}^N}\varphi(y)U_1(|y|)\frac{y}{|y|}\,dy\biggr)\\
 & =\frac{1}{c_1^2}\biggr[2^N N\Gamma\biggr(\frac{N}{2}\biggr)\biggr]^{-1}Nq_*^2
  \biggr(\frac{x}{|x|}\cdot\Xi(\varphi)\biggr).
\end{split}
\end{equation*}
Then we have 
\begin{equation}
\label{eq:3.29}
\begin{split}
 & -q_*\frac{N}{2} M_{0, 1}U(|x|)F_0(|x|)+\sum^N_{i=1}M_{1, i}U_1(|x|)Q_{1, i} \left( \frac{x}{|x|} \right)\\
 & =\frac{q_*^2}{2^N\Gamma(N/2)}
\biggr[-\frac{N}{c_*^2}M(\varphi)+\frac{1}{c_1^2}\biggr(\frac{x}{|x|}\cdot\Xi(\varphi)\biggr)U_1(|x|)\biggr].
\end{split}
\end{equation}
Then Theorem~\ref{Theorem:3.6} follows from \eqref{eq:3.26}, \eqref{eq:3.27}, \eqref{eq:3.28} and \eqref{eq:3.29}. 
$\Box$\vspace{3pt}
\newline
Modifying Theorem~\ref{Theorem:3.5}, 
we give sufficient conditions for the hot spots to 
consist of only one point and to move along a smooth curve in the case where $A=0$ and $\Pi\not=\emptyset$. 
We remark that $A_1=A+1$ if $A=0$. 
\begin{theorem}
\label{Theorem:3.7}
Let $L_V$ be a nonnegative Schr\"odinger operator under condition~{\rm (V)} with 
$\lambda_1\ge 0$,  $A=0$ and $\Pi\not=\emptyset$. 
Assume \eqref{eq:2.17} and let $u$ be a solution of \eqref{eq:1.1}. 
Let $x_*\in{\bf R}^N$ be such that $|x_*|\in\Pi$ and 
$$
\lim_{t\to\infty}\sup_{x\in H(u(t))}|x-x_*|=0. 
$$
Then there exist a constant $T>0$ and a curve $x=x(t)\in C^1([T,\infty):{\bf R}^N)$ such that 
$H(u(t))=\{x(t)\}$ for $t\ge T$ in the following cases: 
\begin{itemize}
  \item[{\rm (a)}] 
  $|x_*|=0$, $V\in C^\gamma([0,\infty))$ for some $\gamma\in(0,1)$ and 
  $\nabla^2 U(|x|)\le 0$ in a neighborhood of $x=0$; 
  \item[{\rm (b)}] 
  $|x_*|>0$,  $U''(r_*)<0$ and $\Xi(\varphi)\not=0$;
  \item[{\rm (c)}] 
  $|x_*|>0$, $U''(r)\le 0$ in a neighborhood of $r=r_*$, $S''(r_*)<0$ and $\Xi(\varphi)\not=0$.
  \end{itemize}
Here $S=S(r)$ is as Theorem~{\rm\ref{Theorem:3.6}}. 
\end{theorem}
{\bf Proof.}
The proofs in case~(a) with $U''(0)<0$ and case~(b) 
are obtained by the same argument as the proof of Theorem~\ref{Theorem:3.5} in cases~(I) and (III), respectively. 
So it suffices to consider case~(a) with $U''(0)=0$ and case~(c).

Let us consider case~(a) with $U''(0)=0$. Let $\epsilon>0$. 
It follows from \eqref{eq:3.1} that $U_1'(0)=0$ and 
$$
U_1''(r)=-(N-1)r^{-N-1}\int_0^r \tau^NV(\tau)U_1(\tau)\,d\tau+V(r)U_1(r)\to 0
\quad\mbox{as}\quad r\to 0.
$$
These imply that $U_1\in C^2([0,\infty))$ and $U_1''(0)=0$. 
Then, similarly to \eqref{eq:3.20}, we have 
$$
-t^{\frac{N+2A_1}{2}}(\partial_r^2u_{1,1})(x,t)\ge -\epsilon
$$
in a neighborhood of $x=0$ for all sufficiently large $t>0$. 
Then, applying a similar argument as in proof of Theorem~\ref{Theorem:3.5} in case~(II), 
we obtain Theorem~\ref{Theorem:3.7} in case~(a) with $U''(0)=0$.

Let us consider case~(c). Similarly to the proof of Theorem~\ref{Theorem:3.5} in case~(III), 
without loss of generality, we can assume that $\Xi(\varphi)/|\Xi(\varphi)|=(1,0,\dots,0)$ and $x_*=(r_*,0,\dots,0)$ 
and we introduce the coordinate $(r,\theta_2,\dots,\theta_N)$ in a neighborhood of $x=x_*$. 
Then, by \eqref{eq:2.3}, \eqref{eq:2.6} and \eqref{eq:3.20} we have  
\begin{equation}
\label{eq:3.30}
\begin{split}
 & -q_*^{-1}t^{\frac{N}{2}}(\partial_r^2 u_0)(x,t)\\
 & =-[M_{1,0}+o(1)](\partial_r^2 U)(|x|)
+\left[\frac{N}{2}M_{0,1}+o(1)\right]t^{-1}\partial_r^2[UF_0](|x|)+O(t^{-2})\\
 & \ge\left[\frac{N}{2}M_{0,1}+o(1)\right]t^{-1}\partial_r^2[UF_0](|x|)+O(t^{-2}),\\
 & -q_N^{-1}t^{\frac{N}{2}+1}(\partial_r^2 u_{1,1})(x,t)=
-[M_{1,1}+o(1)](\partial_r^2 U_1)(|x|)\theta_i+O(t^{-1}),
\end{split}
\end{equation}
in a neighborhood of $x=x_*$ for all sufficiently large $t>0$, where $i=1,\dots,N$. 
By condition~(c), \eqref{eq:3.29} and \eqref{eq:3.30} we obtain  
\begin{equation*}
\begin{split}
 & -t^{\frac{N}{2}+1}[(\partial_r^2 u_0)(x,t)+(\partial_r^2 u_{1,1})(x,t)]\\
 & \ge -\frac{q_*^2}{2^N\Gamma(N/2)}S''(r)+o(1)+O(t^{-1})\ge-\frac{q_*^2}{2^{N+1}\Gamma(N/2)}S''(r_*)>0
\end{split}
\end{equation*}
in a neighborhood of $x=x_*$ for all sufficiently large $t>0$. 
Similarly, we have 
\begin{equation}
\label{eq:3.31}
-t^{\frac{N}{2}+1}(\partial_r\partial_{\theta_\alpha}u_{1,1})(x,t)=O(|\theta_\alpha|)+O(t^{-1})
\end{equation}
in a neighborhood of $x=x_*$ for all sufficiently large $t>0$, where $\alpha=2,\dots,N$. 
Furthermore, similarly to the proof of Theorem~\ref{Theorem:3.5} in case~(III), 
we have \eqref{eq:3.18}, \eqref{eq:3.22} and \eqref{eq:3.23}. 
Then, combining \eqref{eq:3.30} and \eqref{eq:3.31}, 
we see that $-(\tilde{\nabla}^2u)(x,t)$ is positive definite in a neighborhood of $x_*=(r_*,0)$ for all sufficiently large $t>0$. 
Thus Theorem~\ref{Theorem:3.7} in case~(c) follows from the implicit function theorem. 
Therefore the proof of Theorem~\ref{Theorem:3.7} is complete.
$\Box$
\begin{remark}
\label{Remark:3.3}
Consider the case of the heat equation under condition~\eqref{eq:2.17}. 
Then $V\equiv 0$, $c_*=1$, $c_1=1$, $U(r)=1$, $U_1(r)=r$,  $F_0(r)=r^2/(2N)$ and $\Pi=[0,\infty)$. 
Since
$$
S'(r)=-rM(\varphi)+|\Xi(\varphi)|=-r\int_{{\bf R}^N}\varphi(y)\,dy+\biggr|\int_{{\bf R}^N}y\varphi(y)\,dy\biggr|,
$$
it follows from Theorem~{\rm\ref{Theorem:3.6}} 
that the hot spots converges to $\int_{{\bf R}^N}y\varphi(y)\,dy/\int_{{\bf R}^N}\varphi(y)\,dy$. 
Furthermore, if $\Xi(\varphi)\not=0$, then, 
by Theorem~{\rm\ref{Theorem:3.7}~(c)} we see that 
the hot spots consist of only one point and move along a smooth curve. 
These coincide with statements~{\rm (H2)} and {\rm (H3)} in Section~{\rm 1}. 
\end{remark}

\noindent
{\bf Acknowledgements.} 
The first author was partially supported 
by the Grant-in-Aid for Scientific Research (A)(No.~15H02058)
from Japan Society for the Promotion of Science. 
The second author was  
supported in part by JSPS KAKENHI (Grant No.~15K04965 and 15H03631).
 

\begin{thebibliography}{xx}

\bibitem{BE} 
J. Bebernes and D. Eberly, 
{\it Mathematical Problems from Combustion Theory}, 
Mathematical Science, {\bf 83}, Springer-Verlag, New York, 1989.

\bibitem{CK} 
I. Chavel and L. Karp, 
Movement of hot spots in Riemannian manifolds, 
J. Analyse Math. {\bf 55} (1990), 271-286.

\bibitem{CK1} 
I. Chavel and L. Karp, 
Large time behavior of the heat kernel: the parabolic $\lambda$-potential alternative, 
Comment. Math. Helv. {\bf 66} (1991), 541--556. 

\bibitem{Dav}
E. B. Davies, 
{\it Heat Kernels and Spectral Theory}, Cambridge Tracts in Math. 92,
Cambridge Univ. Press 1989.

\bibitem{DS} 
E. B. Davies and B. Simon, 
$L^p$ norms of noncritical Schr\"odinger semigroups, 
J. Funct. Anal. {\bf 102} (1991), 95--115.

%

\bibitem{Gri} 
A. Grigor'yan, 
{\it Heat Kernel and Analysis on Manifolds,} 
AMS, Providence, RI,  2009.

\bibitem{I1}
K. Ishige, 
Movement of hot spots on the exterior domain
of a ball under the  Neumann boundary condition, 
J. Differential Equations {\bf 212} (2005), 394-431.

\bibitem{I2}
K. Ishige,
Movement of hot spots on the exterior domain
of a ball under the  Dirichlet boundary condition, 
Adv. Differential Equations {\bf 12} (2007), 1135-1166.

\bibitem{IK02}
K. Ishige and Y. Kabeya,
Large time behaviors of hot spots for the heat equation with a potential, 
J. Differential Equations {\bf 244} (2008), 2934--2962; 	
{\it Corrigendum} in J. Differential Equations {\bf 245} (2008), 2352--2354. 

\bibitem{IK03}
K. Ishige and Y. Kabeya,
Hot spots for the heat equation with a rapidly decaying negative potential, 
Adv. Differential Equations {\bf 14} (2009), 643--662.
	
\bibitem{IK04}
K. Ishige and Y. Kabeya, 
Hot spots for the two dimensional heat equation with a rapidly decaying negative potential, 
Discrete Contin. Dyn. Syst. Ser. S {\bf 4} (2011), 833--849.
	
\bibitem{IK05} 
K. Ishige and Y. Kabeya, 
 $L^p$ norms of nonnegative Schr\"odinger heat semigroup and the large time behavior of hot spots, 
J. Funct. Anal. {\bf 262} (2012), 2695--2733. 
	
\bibitem{IKO}
K. Ishige, Y. Kabeya and E.~M. Ouhabaz, 
The heat kernel of a Schr\"odinger operator with inverse square potential, 
Proc. Lond. Math. Soc. {\bf 115} (2017), 381--410.

\bibitem{IM}
K. Ishige and A. Mukai, 
Large time behavior of solutions of the heat equation with inverse square potential, 
preprint (arXiv:1709.00809).
 
\bibitem{JS} 
S. Jimbo and S. Sakaguchi, 
Movement of hot spots over unbounded domains in ${\bf R}^N$, 
J. Math. Anal. Appl. {\bf 182} (1994), 810--835.

\bibitem{LSU:LSUbook} 
O. A.~Ladyzenskaja, V. A.~Solonnikov, and N. N.~Ural'ceva, 
{\it Linear and Quasilinear Equations of Parabolic Type}, 
Amer. Math. Soc., Providence, RI, 1968.

\bibitem{LS}
V. Liskevich and Z. Sobol, 
Estimates of integral kernels for semigroups associated with second-order elliptic operators with singular coefficients, 
Potential Anal. {\bf 18} (2003), 359--390. 

\bibitem{Marchi} 
C. Marchi, 
The Cauchy problem for the heat equation with a singular potential, 
Differential Integral Equations {\bf 16} (2003), 1065--1081.

\bibitem{MS1} 
P.~D. Milman and Y.~A. Semenov, 
Global heat kernel bounds via desingularizing weights, 
J. Funct. Anal. {\bf 212} (2004), 373--398. 

\bibitem{MT1} 
L. Moschini and A. Tesei, 
Harnack inequality and heat kernel estimates for the Schr\"odinger operator with Hardy potential, 
Rend. Mat. Acc. Lincei {\bf 16} (2005), 171--180.

\bibitem{MT2} 
L. Moschini and A. Tesei, 
Parabolic Harnack inequality for the heat equation with inverse-square potential, 
Forum Math. {\bf 19} (2007), 407--427.

\bibitem{M} 
M. Murata, 
Structure of positive solutions to $(-\Delta+V)u=0$ in ${\bf R}^n$, 
Duke Math. J. {\bf 53} (1986), 869--943.

\bibitem{Ouh}
E.~M. Ouhabaz, 
{\it Analysis of Heat Equations on Domains}, London Math. Soc. Monographs, 31,
Princeton Univ. Press 2005.

\bibitem{P3} 
Y. Pinchover, 
On criticality and ground states of second order elliptic equations, II, 
J. Differential Equations {\bf 87} (1990), 353--364.

\bibitem{P0} 
Y. Pinchover, 
Large time behavior of the heat kernel and the behavior 
of the Green function near criticality for nonsymmetric elliptic operators, 
J. Funct. Anal. {\bf 104} (1992), 54--70.

\bibitem{P1} 
Y. Pinchover, 
On positivity, criticality, and the spectral radius of the shuttle operator for elliptic operators, 
Duke Math. J. {\bf 85} (1996), 431--445. 

\bibitem{P2} 
Y. Pinchover, 
Large time behavior of the heat kernel, 
J. Funct. Anal. {\bf 206} (2004), 191--209.

\bibitem{PZ} 
Y. Pinchover, 
Some aspects of large time behavior of the heat kernel: an overview with perspectives, 
Mathematical Physics, Spectral Theory and Stochastic Analysis (Basel) (M. Demuth and W. Kirsch, eds.), 
Operator Theory: Advances and Applications, vol. 232, Springer Verlag, 2013,  
299--339.

\bibitem{S} 
B. Simon, 
Large time behavior of the $L^p$ norm of Schr\"odinger semigroups, 
J. Funct. Anal. {\bf 40} (1981), 66--83.

\bibitem{V} 
J.~L. Vazquez, 
Domain of existence and blowup for the exponential reaction-diffusion equation, 
Indiana Univ. Math. J. {\bf 48} (1999), 677--709. 

\bibitem{VZ}
J. L. V\'azquez and E. Zuazua, 
The Hardy inequality and the asymptotic behaviour of the
heat equation with an inverse-square potential, 
J. Funct. Anal. {\bf 173} (2000), 103--153. 

\bibitem{Zhang0}
Qi S. Zhang, 
Large time behavior of Schr\"odinger heat kernels and applications, 
Comm. Math. Phys. {\bf 210} (2000), 371--398.

\bibitem{Zhang}
Q. S. Zhang, 
Global bounds of Schr\"odinger heat kernels 
with negative potentials, 
J. Funct. Anal. {\bf 182} (2001), 344--370.

\bibitem{Z}
Ya. B. Zel'dovich, G.~I. Barenblatt, V.~B. Librovich and G.~M. Makhviladze, 
{\it The mathematical theory of combustion and explosions}, 
Consultants Bureau, New York, 1985. 
\end{thebibliography}
\end{document}